\documentclass[11pt,reqno]{amsart}

\usepackage[T1]{fontenc}
\usepackage{amsmath,amssymb,amsthm}
\usepackage{xcolor}
\usepackage[colorlinks=true,linkcolor=blue,citecolor=blue,urlcolor=blue]{hyperref}

\numberwithin{equation}{section}
\usepackage{dsfont}

\theoremstyle{plain}
\newtheorem{theorem}{Theorem}[section]
\newtheorem{proposition}[theorem]{Proposition}
\newtheorem{lemma}[theorem]{Lemma}
\newtheorem{corollary}[theorem]{Corollary}

\theoremstyle{definition}
\newtheorem{definition}[theorem]{Definition}

\theoremstyle{remark}

\newcommand{\C}{\mathds{C}}
\newcommand{\R}{\mathds{R}}
\newcommand{\Z}{\mathds{Z}}
\newcommand{\CP}{\mathds{C}\mathrm{P}}

\newcommand{\ee}{\mathrm{e}}

\makeatletter
\@namedef{subjclassname@2020}{\textup{2020} Mathematics Subject Classification}
\makeatother

\begin{document}

\title[Stable immersions of radial K\"ahler metrics]{Stable projective immersions of K\"ahler metrics with some radial symmetry}

\author{Claudio Arezzo}
\address{The Abdus Salam International Centre for Theoretical Physics, Strada Costiera 11,
34151 Trieste, Italy}
\email{arezzo@ictp.it}

\author{Andrea Loi}
\address{Dipartimento di Matematica e Informatica, Universit\`a di Cagliari, Via Ospedale 72,
09124 Cagliari, Italy}
\email{loi@unica.it}

\author{Giovanni Placini}
\address{Dipartimento di Matematica e Informatica, Universit\`a di Cagliari, Via Ospedale 72,
09124 Cagliari, Italy}
\email{giovanni.placini@unica.it}

\author{Michela Zedda}
\address{Dipartimento di Scienze Matematiche, Fisiche e Informatiche, Universit\`a di Parma,
Parco Area delle Scienze 53/A, 43124 Parma, Italy}
\email{michela.zedda@unipr.it}



\begin{abstract}
We combine Calabi’s criterion with Pringsheim’s theorem to study projective inducibility under positive rescaling. By restricting to suitable radial holomorphic curves, we obtain an obstruction to local K\"ahler immersions into finite- or infinite-dimensional complex projective space that applies uniformly to every positive rescaling. For invariant rational curves, we show that non-rationality of the moment map implies that no positive rescaling is projectively induced. We apply these criteria to the Taub–NUT, Candelas–de la Ossa, non-flat Stenzel, and Hwang–Singer scalar-flat metrics, as well as to Gibbons--Hawking metrics.
Except for the Burns--Simanca metric in the Hwang--Singer family, no positive rescaling is projectively induced. For Gibbons--Hawking ALE spaces with at least two centres, the conclusion holds in every compatible complex structure.  We also extend the results to complex space forms
of arbitrary curvature with at least one finite signature index.
\end{abstract}

\maketitle

\tableofcontents

\section{Introduction and statement of the main results}

A K\"ahler metric $g$ on a complex manifold $M$ is called {\em projectively induced} if it is induced by a local K\"ahler immersion into some $\CP^N$ endowed with its Fubini--Study metric, with $N\in\Z_{\geq0}\cup\{\infty\}$.

A basic question in the theory of projectively induced K\"ahler metrics is how projective
inducibility behaves under rescaling. Given a K\"ahler metric $\omega$, we ask
for which $c>0$ the metric $c\,\omega$ is projectively induced. In the compact case when the immersion is global the admissible scales are necessarily discrete. A global projective immersion forces the K\"ahler class of $c\,\omega$ to be integral, up to the fixed normalization of the Fubini--Study form,
so the ray $\{c[\omega]:c>0\}$ can meet the integral lattice only discretely. The basic example
is the Fubini--Study metric on $\CP^n$, for which $c\,\omega_{FS}$ is projectively induced
precisely when $c\in\mathbb Z_{>0}$, through the Veronese embeddings.

In the non-compact case this cohomological quantization could disappear and the dependence on the
scale can be much less rigid. 
For example every positive multiple of the flat metric on $\C^n$ is projectively induced into $\CP^\infty$. This example shows that an entire ray of positive multiples can be projectively induced for a noncompact metric. In
principle, each multiple of a K\"ahler metric could behave differently. For example, for an irreducible bounded symmetric domain endowed with its Bergman metric, the admissible rescalings are governed by the Wallach set (see \cite{loi-zedda-wallach} and also \cite{DSIL} for the homogeneous case), and --in rank greater
than one-- this set has both a finite discrete part and a continuous half-line. This contrast
makes it natural to seek criteria that determine the whole set of admissible scales without
treating each value of $c$ separately.

Calabi's resolvability criterion~\cite{calabi} gives an algebraic characterization of local
projective inducibility in terms of the positivity of an infinite coefficient matrix. In
principle it can therefore be used to test each rescaled metric $c\,\omega$ separately. For
the rescaling problem, however, this direct approach rapidly becomes impractical. One must
find a negative coefficient or minor for each value of $c$, and in concrete non-compact
examples the first obstruction may occur at progressively higher order as the scale varies.
Thus a coefficient-by-coefficient use of Calabi's criterion gives no effective uniform
control over the whole range $c>0$. This phenomenon has been observed for example for the Taub--NUT metric in \cite{lzz-taubnut}.

The first main result of this paper gives a uniform obstruction designed precisely to overcome this difficulty
and treat all positive multiples at once. The idea is to restrict the metric to a holomorphic
curve on which the induced metric is radial (\emph{radial holomorphic curves} in the following). Along such a curve, the relevant entries of
Calabi's matrix become the coefficients of a one-variable power series. If a positive
multiple of the metric is projectively induced, these coefficients must be non-negative.
Pringsheim's Theorem then forces a dominant singularity on the positive real axis. A
logarithmic-derivative argument makes this singularity obstruction independent of the scale.

\begin{theorem}\label{thm:main1}
Let $(M,\omega)$ be a real-analytic K\"ahler manifold admitting a radial holomorphic curve
with potential $\Phi$. If $\Phi$ extends analytically along the positive real axis
whereas $\Phi'$ has no single-valued meromorphic continuation to $\C$, then no positive multiple of $\omega$ admits a local K\"ahler immersion into $\CP^N$, for any $N\leq\infty$.
\end{theorem}

We apply Theorem \ref{thm:main1} to families of manifolds with significant curvature properties.
The case of Ricci-flat K\"ahler metrics is of particular interest. In \cite{Hulin}, Hulin proved that compact K\"ahler--Einstein submanifolds of $\CP^N$, for finite $N$, have strictly positive Einstein constant. Notice that for compact manifolds restricting to finite dimensional ambient space is not reductive, since a compact K\"ahler submanifold of $\CP^\infty$ is actually contained in some $\CP^N$ of finite dimension.  In contrast, the non-compact case is more flexible and we have both examples of K\"ahler submanifolds of finite dimensional complex projective spaces, and K\"ahler submanifolds of $\CP^\infty$ that are {\em full}, in the sense that their image is not contained in any finite dimensional $\CP^N$. In \cite{ALL}, Arezzo, Li and Loi exclude global K\"ahler immersions of complete noncompact Ricci-flat manifolds into finite-dimensional projective space. 
The conjecture that every projectively induced Ricci-flat K\"ahler metric is flat appeared in \cite{lsz}, where the authors proved it for radial K\"ahler metrics under a stability assumption. Subsequently the stability assumption has been removed in \cite{ALPZradial}, and it has been verified for several other classes, including Calabi's Ricci-flat metrics \cite{lzz-calabi} and K\"ahler cones over regular complete Sasakian manifolds \cite{MTZ}. In general, the conjecture remains open.
The scalar-flat case has a different feature, as shown for example by the Burns--Simanca metric, which is proven to be projectively induced in \cite{cannas-loi,lsz}. 

These two curvature classes therefore provide complementary settings in which to test the uniform obstruction.

\begin{theorem}\label{thm:main4}
For each of the following metrics, no positive multiple admits a local K\"ahler immersion into $\CP^N$ for any $N\leq\infty$:
\begin{enumerate}
\item[(i)] the Taub--NUT metric $\omega_m$ on $\C^2$, for every $m>0$;
\item[(ii)] the Candelas--de la Ossa Ricci-flat metrics;
\item[(iii)] the Stenzel Ricci-flat non-flat metrics on the complexifications of compact rank-one symmetric spaces;
\item[(iv)] the Hwang--Singer scalar-flat ALE metric on $\mathcal O_{\CP^{n-1}}(-k)$, for every $n\ge2$ and $k\ge1$ with $(n,k)\neq(2,1)$;
\item[(v)] the Hwang--Singer scalar-flat metric on $\C^{n+1}$ associated with $\beta < 0$;
\item[(vi)] the Gibbons--Hawking ALE metric determined by at least two distinct centres, in any compatible complex structure.
\end{enumerate}
\end{theorem}

Cases $(i)$, $(ii)$, $(iii)$, and $(vi)$, provide further evidence for the Ricci-flat conjecture. We stress that our results concern the entire rescaling ray, indeed rather than testing projective inducibility at a fixed scale, they exclude all positive multiples simultaneously. For the Taub--NUT metrics, Loi--Zedda--Zuddas excluded projective immersions when $m>c/2$ \cite{lzz-taubnut}; the all-scale Taub--NUT problem has been recently solved by Huang in an AI generated preprint \cite{HuangTaubNut}. Accordingly, the all-scale Taub--NUT conclusion is included here as an application of our common criterion, without a claim of exclusive novelty. To the best of our knowledge, the Candelas--de la Ossa metrics had not previously been studied in the projective-rescaling problem. Stenzel metrics were investigated in \cite{zedda-stenzel}, where the metrics on the complexifications of $\CP^n$ and $\mathds{H}{\rm P}^n$ were shown not to be projectively induced when rescaled by a factor $0<c\leq 1$. Our result gives a uniform treatment of the Stenzel metrics on the complexifications of all compact rank-one symmetric spaces and excludes every positive rescaling.

The scalar-flat ALE metrics obtained from the Hwang--Singer momentum construction~\cite{hwang-singer} are $U(n)$-invariant and exhibit a different phenomenon. Within this family one recovers the Burns--Simanca metric, which is projectively induced only if rescaled by an integral factor \cite{cannas-loi}, and in \cite{cristofori-zedda} it was characterized within the Hwang--Singer family on $\mathcal O_{\CP^1}(-k)$ as the only projectively induced metric. For the metrics on $\C^{n+1}$, Cristofori--Zedda already excluded projective immersion of \emph{every} positive multiple \cite[Theorem~1]{cristofori-zedda}; case (v) gives another proof within the present framework. For the line-bundle family they characterized projective inducibility at the original scale \cite[Theorem~2]{cristofori-zedda}, while all-scale non-inducibility of the Eguchi--Hanson member was known from \cite{lzz-calabi}. Our new projective contribution for the other line-bundle members is the exclusion of every positive multiple. Thus, within this natural family, the Burns--Simanca metric emerges as the unique exception to an otherwise uniform non-inducibility phenomenon.

The Gibbons--Hawking examples fit the same one-variable obstruction. For every compatible complex structure, the circle-fibration preimage of an outward ray from a suitably chosen centre defines a radial holomorphic curve. When there are at least two centres, the induced potential satisfies the hypotheses of Theorem \ref{thm:main1}, yielding non-inducibility for every positive rescaling and every compatible complex structure. In the two-centre case, this includes the known non-inducibility of the Eguchi--Hanson metric in its resolution complex structure \cite{lzz-calabi}. The construction and proof are given in Section 6. A further consequence of Theorem \ref{thm:main1} is the following obstruction in terms of the moment map of an invariant rational curve.

\begin{theorem}\label{thm3}
Let $(M,\omega)$ be a real-analytic K\"ahler manifold containing a holomorphic curve $\Sigma \simeq \mathbb{CP}^{1}$
whose induced metric is invariant under a non-trivial and holomorphic $S^{1}$-action. If the moment map of this action is non-rational as a function of $\rho=|z|^2$ in an affine holomorphic coordinate $z$ adapted to the action, then, for every $c>0$, the K\"ahler manifold $(M,c\omega)$ admits no local K\"ahler immersion into $\CP^N$ for any $N\leq\infty$.
\end{theorem}

 Although we focus on projective immersions, the same method also
excludes local K\"ahler immersions into complex space forms of
arbitrary curvature whenever at least one signature index is finite. 
This restriction is sharp, since every real-analytic K\"ahler metric
admits local immersions into such spaces when both indices are infinite. We treat the indefinite case in the last section.\\

The paper is organized as follows. Section~\ref{sec:calabi} recalls Calabi's projective criterion and its one-variable reduction. Section~\ref{sec:sharp} establishes the Pringsheim obstruction and proves Theorem~\ref{thm:main1}. Section~\ref{sec:applications} treats the first five families in Theorem~\ref{thm:main4}. Section~\ref{sec:moment} develops the invariant-sphere criterion, and Section~\ref{sec:GH} applies Theorem~\ref{thm:main1} to Gibbons--Hawking metrics and all their compatible complex structures. Finally, in Section~\ref{7} we extend our results to complex space forms of arbitrary curvature and signature.

\section{Calabi's criterion on a holomorphic curve}\label{sec:calabi}

We recall only the part of Calabi's theory needed below. For this part we refer the reader to \cite{calabi} and to
\cite[Ch.~2]{lz-book} for details and further results. Let $(M,\omega)$ be a real-analytic K\"ahler manifold and
let $p\in M$. If $\Phi$ is a K\"ahler potential near $p$ and $\widetilde\Phi(z,\bar u)$
denotes its analytic continuation to a neighbourhood of the diagonal of $M\times \bar M$, the diastasis of
$\omega$ centred at $p$ is
\begin{equation}\label{eq:diastasis}
D_p(z,\bar z)=\widetilde\Phi(z,\bar z)+\widetilde\Phi(p,\bar p)
-\widetilde\Phi(z,\bar p)-\widetilde\Phi(p,\bar z),
\end{equation}
and in local coordinates centred at $p$ it has an expansion
\[
D_p(z,\bar z)=\sum_{\alpha,\beta}d_{\alpha\beta}\,z^\alpha\bar z^\beta
\]
with no purely holomorphic or anti-holomorphic terms.

In an affine chart centred at the origin, the Fubini--Study diastasis on $\CP^N$ is $D_{FS}(w)=\log(1+\sum_{j=1}^{N}|w_j|^2)$. Thus the matrix relevant to projective immersion is the coefficient matrix of $e^{D_p}-1$, as shown by Calabi's criterion:

\begin{theorem}[Calabi \cite{calabi}]\label{thm:calabi}
A real-analytic K\"ahler manifold $(M,\omega)$ admits a local K\"ahler immersion of a neighborhood of $p\in M$  into $\CP^N$, $N\leq\infty$, if and only if the Hermitian matrix $\bigl(c_{\alpha\beta}\bigr)$ defined by
\begin{equation}\label{defmatrixdiast}
ìe^{D_p}-1=\sum_{\alpha,\beta}c_{\alpha\beta}z^\alpha\bar z^\beta
\end{equation}
is positive semidefinite of rank at most $N$. 
\end{theorem}
\begin{definition}\label{def:radialcurve}
We say that $(M,\omega)$ admits a \emph{radial holomorphic curve with potential $\Phi$} if there exist an open set
$\Omega\subseteq\C$ with $0\in\Omega$, a holomorphic immersion
$\gamma\colon\Omega\to M$, and a real-analytic real-valued function $\Phi$ defined on an interval containing $0$, such that
$\Phi(|z|^2)$ is a K\"ahler potential for $\gamma^*\omega$ near $z=0$. 
\end{definition}

The following proposition guarantees the existence of a radial holomorphic curve in a special case.
\begin{proposition}
Let $(M,\omega)$ be a real-analytic K\"ahler manifold carrying a holomorphic isometric $S^1$-action and let $p\in M$ be a point fixed by this action. Assume that the isotropy representation on $T_pM$ has a non-zero weight. Then $M$ admits a radial holomorphic curve.  
\end{proposition}
\begin{proof}
Let $\phi$ be a real-analytic K\"ahler potential for $\omega$ around $p$, defined on a $S^1$-invariant domain (to obtain a $S^1$-invariant domain take a sufficiently small geodesic ball centered at $p$ contained in the domain of $\phi$). For $e^{i\theta}\in S^1$ let $a_\theta$ denote the action. Since $S^1$ acts by holomorphic isometries, $a_\theta^*\omega=\omega$ and $\phi\circ a_\theta$ is also a K\"ahler potential for $\omega$ around $p$. Taking the average:
$$
\varphi(q):=\frac1{2\pi}\int_{0}^{2\pi}\phi(a_\theta q)d\theta
$$
we obtain a K\"ahler potential which is also $S^1$ invariant:
$$
\frac i2\partial\bar\partial\frac1{2\pi}\int_{0}^{2\pi}\phi(a_\theta q)d\theta=\frac1{2\pi}\int_{0}^{2\pi}a_\theta^*\left(\frac i2\partial\bar\partial \phi(q)\right)d\theta=\omega,
$$
$$
\varphi(a_\beta q)=\frac1{2\pi}\int_{0}^{2\pi}\phi(a_\theta a_\beta q)d\theta=\varphi(q).
$$

By the holomorphic Bochner linearization theorem
(see, for instance, \cite[Lemma~2.2]{TW}), there exist neighborhoods \(U\) of \(p\) and \(V\) of \(0\in T_pM\),
and an \(S^1\)-equivariant biholomorphism
\[
F:U\longrightarrow V
\]
with \(F(p)=0\), such that the \(S^1\)-action on \(V\) coincides with
the isotropy representation. After shrinking $U$ if necessary, we may
assume that $U$ is contained in the domain of $\varphi$.
 In coordinates adapted to the
weight decomposition, the action takes the form
\[
e^{i\theta}\cdot(z_1,\ldots,z_n)
=
\bigl(
e^{im_1\theta}z_1,\ldots,e^{im_n\theta}z_n
\bigr),
\]
where \(m_1,\ldots,m_n\in\mathbb Z\).
Since the isotropy representation has a nonzero weight, there exists
an index \(j\) such that \(m_j\neq 0\). Then, denoting $e_j$ the $j$-th coordinate vector,
\[
e^{i\theta}\cdot e_j=e^{im_j\theta}e_j.
\]
Set $\Delta_\epsilon:=\{w\in \C|\, |w|<\epsilon\}$. For small enough $\epsilon$, consider the holomorphic curve $\gamma\!:\Delta_\epsilon\rightarrow M$, given by
$$
\gamma(w)=F^{-1}(we_j).
$$
Observe that $\gamma$ is a holomorphic immersion, in fact, $\gamma'(w)=d(F^{-1})_{we_j}(e_j)\neq 0$, throughout the chosen disc, because $F^{-1}$ is a biholomorphism, and $\gamma(0)=p$. 
Then
$$
e^{i\theta}\cdot\gamma(w)= \gamma(e^{im_j\theta}w).
$$
Since $\varphi$ is $S^1$-invariant,
$$
\varphi(\gamma(w))= \varphi(e^{i\theta}\cdot\gamma(w)) =\varphi(\gamma(e^{im_j\theta}w)).
$$
Since $m_j\neq 0$, the map $S^1\rightarrow S^1$, $e^{i\theta}\mapsto e^{im_j\theta}$, is surjective. Therefore
$$
\varphi(\gamma(e^{i\beta}w))=\varphi(\gamma(w))
$$
for every $\beta\in\R$, that is $\varphi\circ\gamma$ is rotation invariant.

Since $\varphi\circ\gamma$ is real analytic, it admits an expansion
$$
\varphi(\gamma(w))=\sum_{r,s\geq 0} c_{rs}w^r\bar w^s,
$$
and rotation invariance gives
$$
\sum_{r,s\geq 0}
c_{rs}e^{i(r-s)\beta}w^r\bar w^s=\sum_{r,s\geq 0}c_{rs}w^r\bar w^s
$$
for every $\beta$. Hence $c_{rs}=0$ whenever $r\neq s$. It follows that
$$
\varphi(\gamma(w))=\sum_{r\geq 0} c_{rr}|w|^{2r}=\Phi(|w|^2)
$$
for some real-analytic function $\Phi$. Thus, $\gamma$ is a radial holomorphic curve with potential $\Phi$, as wished.
\end{proof}

For a radial potential, Calabi’s matrix is diagonal. Restricting to a radial holomorphic curve therefore gives the following necessary condition.

\begin{lemma}\label{lem:diagonal}
Assume that $(M,\omega)$ admits a radial holomorphic curve with potential $\Phi$. If $\omega$ admits a local K\"ahler immersion into $\CP^N$ and $e^{\Phi(x)}=\sum_{j\geq0}a_jx^j$, then $a_j\geq0$ for all $j\geq0$. 
\end{lemma}
\begin{proof}
If $\omega$ is projectively induced then so is $\gamma^*\omega$, and Calabi's criterion (Theorem \ref{thm:calabi} above) implies that the matrix of coefficients of $e^{D(x)}-1$ is positive semidefinite, where $D(x)$ is the diastasis of the induced metric on the curve centred at the origin. Since $D(x)=\Phi(x)-\Phi(0)$, we have $e^{D(x)}-1=e^{-\Phi(0)}e^{\Phi(x)}-1$ and the positive-degree coefficients of $e^D-1$ are positive multiples of those of $e^\Phi$. Since the matrix of coefficients of $e^D-1$ is diagonal, its positive semidefiniteness implies $a_j\geq 0$ for any $j\geq 1$ ($a_0=e^{\Phi(0)}>0$). 
\end{proof}

\section{Pringsheim's Theorem and K\"ahler immersions}\label{sec:sharp}

For a holomorphic germ $f$ at the origin, we denote by
$R(f)\in(0,+\infty]$ the radius of convergence of its Taylor series and by
$$
\rho_+(f)=\sup\{r>0: f \text{ continues analytically along }[0,r]\}.
$$
 Since $f$ is
holomorphic on the disc of radius $R(f)$, one always has
\begin{equation}\label{eq:trivial}
R(f)\le\rho_+(f).
\end{equation}
The following proposition is a consequence of the classical Pringsheim's Theorem (\cite{flajolet-sedgewick}, Theorem IV.6), which states that a power series with non-negative coefficients and finite  radius of convergence $R<\infty$, has a singularity at $R$. The same conclusion holds if the coefficients are eventually
non-negative, since adding a polynomial changes neither the
radius of convergence nor the singularities of the germ.
\begin{proposition}\label{cor:Rrho-calabi}
Let $(M,\omega)$ be a real-analytic K\"ahler manifold admitting a radial holomorphic curve with potential $\Phi$. If $\omega$ admits a local K\"ahler immersion into $\CP^N$ for some $N\leq\infty$, then
$$
R\bigl(e^{\Phi}\bigr)=\rho_+\bigl(e^{\Phi}\bigr).
$$
\end{proposition}

\begin{proof}
By \eqref{eq:trivial} it suffices to prove $\rho_+(e^\Phi)\leq R(e^\Phi)$. There is nothing to prove if the radius is infinite. Otherwise, Lemma~\ref{lem:diagonal} shows that all positive-degree Taylor coefficients of $e^\Phi$ are non-negative. Pringsheim's Theorem makes $x=R(e^\Phi)$ a singularity, and the desired inequality follows.
\end{proof}

In other words, Proposition~\ref{cor:Rrho-calabi} states that a dominant singularity cannot lie strictly closer to the origin than the first obstruction on the positive real axis, as formalized by the following corollary:

\begin{corollary}\label{cor:closer}
Let $(M,\omega)$ be a real-analytic K\"ahler manifold admitting a radial holomorphic curve with potential $\Phi$. Suppose that the germ $e^{\Phi}$ admits analytic continuation along the straight segment
$[0,z_0)$ but not through its endpoint $z_0$, where $z_0\in \C\setminus \R_{\geq0}$ and
$$
|z_0|<\rho_+(e^{\Phi}).
$$
Then $\omega$ admits no local K\"ahler immersion into $\CP^N$ for any $N\leq\infty$.
\end{corollary}
\begin{proof}
Indeed, $R(e^\Phi)\leq |z_0|<\rho_+(e^\Phi)$, contradicting Proposition~\ref{cor:Rrho-calabi}.
\end{proof}

The next corollary gives a useful consequence of Proposition \ref{cor:Rrho-calabi}.

\begin{corollary}\label{cor:nonentire}
Let $(M,\omega)$ be a real-analytic K\"ahler manifold admitting a radial holomorphic curve with potential $\Phi$. Suppose that $\Phi$ continues analytically along $[0,+\infty)$ and that $e^{\Phi}$ is not the germ of an entire function. Then $\omega$ admits no local K\"ahler immersion into $\CP^N$ for any $N\leq\infty$.
\end{corollary}

\begin{proof}
Since $e^\Phi$ continues wherever $\Phi$ does, $\rho_+(e^\Phi)=+\infty$, whereas non-entireness gives $R(e^\Phi)<\infty$. This contradicts Proposition~\ref{cor:Rrho-calabi}.
\end{proof}

We can now prove Theorem \ref{thm:main1}.

\begin{proof}[Proof of Theorem~\ref{thm:main1}]

Fix $c>0$ and suppose that $c\omega$ is projectively induced. Since $\Phi$ continues analytically along the positive real axis, Proposition \ref{cor:Rrho-calabi} applied to $c\omega$ gives:
$$
R\left(e^{c\Phi}\right)=\rho_+\left(e^{c\Phi}\right)=+\infty.
$$
Hence, the germ $e^{c\Phi}$ extends to an entire function. Since
$$
\Phi'=\frac1c\frac{\left(e^{c\Phi}\right)'}{e^{c\Phi}}
$$
near the origin, its logarithmic derivative provides a single-valued meromorphic continuation of $\Phi'$ to $\C$, contradicting the hypothesis.

\end{proof}

\section{Applications}\label{sec:applications}

\subsection{Taub--NUT metrics}\label{sub:taubnut}
The Taub--NUT metrics introduced by LeBrun in \cite{lebrun} form a one-parameter family $\omega_m$, $m\geq0$, of
complete Ricci-flat K\"ahler metrics on $\C^2$. In the complex coordinates
$(z_1,z_2)$ they admit the K\"ahler potential
\[
\Phi_m(u,v)=u^2+v^2+m(u^4+v^4),
\]
where the functions $u$ and $v$ are implicitly determined by
\[
|z_1|=\ee^{m(u^2-v^2)}u,
\qquad
|z_2|=\ee^{m(v^2-u^2)}v.
\]
For $m=0$ one recovers the flat metric, while for $m>0$ these metrics are Ricci--flat
non-flat. In \cite{lzz-taubnut} it is proven that $c\omega_m$ is not projectively induced for any $c<2m$. 

\begin{theorem}\label{thm:TN}
For every $m>0$, no positive multiple of the Taub--NUT metric $\omega_m$ admits a local K\"ahler immersion into $\CP^N$, for any $N\leq\infty$.
\end{theorem}

\begin{proof}
We apply Theorem \ref{thm:main1} using the properties of Lambert $W$ function.
Consider the holomorphic curve $z\mapsto(z,0)$, where we have $v=0$, $|z|^2=u^2e^{2mu^2}$ and $\Phi_m(u,v)=u^2+mu^4$. Thus, setting $x=|z|^2$, the radial potential of the induced metric on the curve at $(0,0)$ reads:
\begin{equation}\label{eq:TNdiast}
\Phi_m(x)=\frac1{2m}W_0(2mx)+\frac1{4m}W_0(2mx)^2 ,
\end{equation}
where $W_0$ denotes the principal branch of the Lambert $W$ function defined by $s=W(s)e^{W(s)}$.
 Differentiating and using that $W_0'(s)=W_0(s)/\bigl(s(1+W_0(s))\bigr)$, gives
\begin{equation}\label{eq:TNder}
\Phi_m'(x)=\frac{W_0(2mx)}{2mx}.
\end{equation}
The principal branch $W_0$ is holomorphic in a neighborhood of the non-negative real half-line, so \eqref{eq:TNdiast} extends analytically along the positive real axis. 
Furthermore, $W_0$ has a square-root branch point at
$-\ee^{-1}$, hence it does not admit a
single-valued meromorphic continuation to $\C$.
If $\Phi'_m(x)$ had such a continuation, multiplication by $2mx$ would give one for $W_0(2mx)$, which is impossible. Conclusion follows by Theorem \ref{thm:main1}.
\end{proof}

\subsection{The Candelas--de la Ossa metric}\label{sub:CdO}
In \cite{cdlo} Candelas and de la Ossa construct a one-parameter family $\omega_{a}$, $a>0$, of complete Ricci-flat
K\"ahler metrics on the small resolution of the conifold
$$
\mathcal{C} =\left\{ w\in\C^{4}\mid w_{1}^{2}+w_{2}^{2}+w_{3}^{2}+w_{4}^{2}=0 \right\}.
$$

We summarize here Candelas and de la Ossa's construction. Consider the matrix
\[
W
=
\frac{1}{\sqrt2}
\begin{pmatrix}
w^3+iw^4 & w^1-iw^2\\
w^1+iw^2 & -w^3+iw^4
\end{pmatrix}
=
\begin{pmatrix}
X&U\\
V&Y
\end{pmatrix},
\]
then the conifold is given by
\[
\det W=XY-UV=0.
\]

One of the two small resolutions $\nu\!: \tilde{\mathcal C}\rightarrow \mathcal C$ is obtained by
replacing the equation $XY-UV=0$ by
\[
\begin{pmatrix}
X&U\\
V&Y
\end{pmatrix}
\begin{pmatrix}
\lambda_1\\
\lambda_2
\end{pmatrix}
=0,
\quad [\lambda_1:\lambda_2]\in\CP^1.
\]
The resulting complex manifold is
\[
\widetilde{\mathcal C}
\cong
\mathcal O_{\CP^1}(-1)\oplus\mathcal O_{\CP^1}(-1),
\]
and the exceptional set
\[
E=\nu^{-1}(0)
\]
is isomorphic to $\CP^1$.

On the affine chart $\{\lambda_1\neq0\}$ set $\lambda={\lambda_2}/{\lambda_1}$.
The equations of the resolution give
\[
X=-\lambda U,
\qquad
V=-\lambda Y,
\]
so $(U,Y;\lambda)$ are holomorphic coordinates on this chart.
Set
\[
\sigma=|U|^2+|Y|^2,
\qquad
\Lambda=1+|\lambda|^2,
\]
and let
\begin{equation}\label{eq:CdlO-radial}
r^2=\operatorname{tr}(W^* W)=\sigma\Lambda.
\end{equation}

The parameter $a$ is the radius of the exceptional $\CP^{1}$.
For $a>0$, the metric $\omega_a$ can be described in a local chart by the K\"ahler potential
\begin{equation}\label{eq:CdlO-pot}
\mathcal K_a
=
\mathcal F_a(r^2)+4a^2\log\Lambda,
\end{equation}
where radial function
$\mathcal F_a$ is determined by
\begin{equation}\label{eq:CdlO}
y^3+6a^2y^2=r^4,\qquad y=r^2\mathcal F_a'(r^2).
\end{equation}
We always take the branch satisfying $y>0$ for $r^2>0$.

\begin{theorem}\label{thm:CdO}
For every $a>0$, no positive multiple of the Candelas--de la Ossa
Ricci-flat K\"ahler metric admits a local K\"ahler immersion into $\CP^N$, for any $N\leq\infty$.
\end{theorem}

\begin{proof}
In order to apply Theorem \ref{thm:main1}, let us choose a holomorphic radial curve. Fix $r_0>0$ and consider the points
$p:=\frac{r_0}{\sqrt2}(1,i,0,0)$ and $v:=\frac{1}{\sqrt2}(0,0,1,i)\in \mathcal C$. Define the curve $\gamma\!:\C\rightarrow \mathcal C$ by $\gamma(\zeta)=p+\zeta v$. Observe that 
$$
|\gamma(\zeta)|^2=r_0^2+|\zeta|^2>0,
$$
so $\gamma$ never meets the singular point of the conifold. 

In terms of the matrix $W$, the curve is
\[
W(\gamma(\zeta))
=
\begin{pmatrix}
0&r_0\\
0&-\zeta
\end{pmatrix},
\]
thus its lift to the small resolution is given on the chart $\{\lambda_1\neq0\}$ by
\[
U=r_0,
\qquad
Y=-\zeta,
\qquad
\lambda=0.
\]
Consequently on $\gamma$
\[
\Lambda=1,
\qquad
\sigma=r_0^2+|\zeta|^2,
\qquad
r^2=r_0^2+|\zeta|^2,
\]
and we obtain:
\[
\gamma^*\mathcal K_a
=
\mathcal F_a(r_0^2+|\zeta|^2).
\]

Setting $x:=|\zeta|^2$,  the radial potential of the restricted metric is therefore
\begin{equation}\label{eq:CdlO-diast}
\Phi_{a,r_0}(x) =\mathcal F_a\bigl(r_0^2+x\bigr).
\end{equation}

To apply Theorem \ref{thm:main1} we need to show that $\Phi_{a,r_0}$ extends analytically along the positive real axis, while $\Phi_{a,r_0}'$ has no single-valued meromorphic continuation to $\C$.

Set $s:=r^2$. From \eqref{eq:CdlO} define
$$
P(y,s):=y^2(y+6a^2)-s^2
$$
and observe that
$$
\frac{\partial}{\partial y}P(y,s)=3y^2+12a^2y
$$
which is positive in the branch where $y(s)>0$ for $s>0$. It follows by the implicit function theorem that $y(s)$ extends holomorphically along the positive real half-line. By \eqref{eq:CdlO} we get that $\Phi_{a,r_0}'$, and thus also $\Phi_{a,r_0}$, extends holomorphically along the positive real half-line.

It remains to show that $\mathcal F_a'(s)$ has no single-valued meromorphic continuation to $\C$. 
Write
\[
s=\psi(y):=y\sqrt{y+6a^2},
\]
where the square root is positive for $y>-6a^2$.
We have
$$
\psi'(y)=\frac{3(y+4a^2)}{2\sqrt{y+6a^2}},\qquad \psi''(y)=\frac34\frac{y+8a^2}{(y+6a^2)^{3/2}}.
$$
Thus the geometric branch $y(s)>0$ for $s>0$ extends through $s=0$, since $\psi'(0)=\sqrt6\,a\neq0$. Continuing it to negative values of $s$, one
reaches the critical point
\[
y_0=-4a^2,
\qquad
s_0=\psi(y_0)=-4\sqrt2\,a^3.
\]
Since
$$
\psi''(y_0)=\frac{3}{2\sqrt2\,a}\neq0,
$$
a Taylor expansion of $\psi(y)$ around $y_0$ gives:
$$
\psi(y)-\psi(y_0)=\frac{1}{2}\psi''(y_0)(y-y_0)^2+O((y-y_0)^3),
$$
hence the inverse branch has a square-root singularity at $s=s_0$,
\begin{equation}\label{cdlofo}
y(s)-y_0=C(s-s_0)^{1/2}+O(s-s_0),
\qquad C\neq0.
\end{equation}
Since $s_0\neq0$, the same square-root term occurs in
\begin{equation}\label{faprimo}
\mathcal F_a'(s)=\frac{y(s)}s.
\end{equation}
Thus $\mathcal F_a'$ has non-trivial monodromy around $s_0$ and cannot admit a single-valued meromorphic continuation to $\C$. Consequently $\Phi_{a,r_0}'(x)$
has the corresponding branch point at $x_0=-r_0^2-4\sqrt2\,a^3$. Conclusion follows by Theorem \ref{thm:main1}.
\end{proof}

\subsection{Stenzel metrics}
\label{sub:stenzel}

Let $M$ be a compact rank-one symmetric space and let $M^{\C}$ denote its
complexification. Stenzel~\cite{stenzel} constructed a complete Ricci-flat K\"ahler
metric on $M^{\C}$ whose K\"ahler potential depends only on the rank-one radial
variable. 

We describe here what we need about the simply connected cases, that is
$S^d$, $d\geq 2$, $\CP^n$, $\mathds{H}{\rm P}^n$, and $\mathbb{O}{\rm P}^2$ (the one-dimensional spherical case is excluded, since the corresponding Stenzel solution is flat).
Let $d=\dim_{\R}M$. Stenzel's metrics are described locally by a K\"ahler potential $f(\mathcal N)$, where $\mathcal N$ is a function of the coordinates that depends on $M$ and $f$ satisfies the ODE:
$$
\mathcal N^{p/2}(f')^{d}(2\mathcal N-1)+2\mathcal N^{1+p/2}(f')^{d-1}(\mathcal N-1)f''=1.
$$
Let $p+q=d-1$. For the projective rank-one symmetric spaces,
$p$ and $q$ are the multiplicities $m_1$ and $m_2$ of the simple
restricted root and its double, respectively. Thus
\[
(p,q)=(2n-2,1),\qquad (4n-4,3),\qquad (8,7)
\]
for $\CP^n$, $\mathds H{\rm P}^n$, and $\mathds O{\rm P}^2$, respectively.
In the spherical case the standard multiplicities are $(d-1,0)$.
After replacing the standard radial variable $T$ by $T=2t$, the
equation takes the same form with $(p,q)=(0,d-1)$.

Set $N=\cosh^2 t$ and $\rho(t):=f(\cosh^2t)$ so that the ODE can be written as
\begin{equation}\label{eq:stenzel-general-ode}
\rho''(t)(\rho'(t))^{d-1}=C\,\sinh^{p}t\,\sinh^{d-1-p}(2t),\qquad C>0,
\end{equation}
(cfr. \cite[Example~2.4]{BD}).

 \begin{lemma}\label{curve}
 Each simply connected complexification of a compact rank-one symmetric space with its Stenzel metric admits a radial holomorphic curve with potential $\Phi(|z|^2)=f(1+|z|^2)$. 
 \end{lemma}
\begin{proof}
We deal with each case separately.
\begin{enumerate}
\item $M=S^d$, $d\geq 2$.
The complexification is the affine quadric
$$
Q_d=\left\{Z\in\C^{d+1}:Z_0^2+\cdots +Z^2_{d}=1\right\}.
$$
Set $\mathcal N_Q(Z)=\sum_{j=0}^d |Z_j|^2$.
Then, with the normalization used above,
\[
N=\frac{\mathcal N_Q+1}{2}.
\]
Hence the radial holomorphic curve is $\gamma(z)=(1,z,iz,0,\ldots,0)$.
\item $M=\CP^n$. In the affine coordinates $(z,w)$ of the complexification used by Lee~\cite{Lee}, the invariant function is
\[
N(z,w)
=
\frac{(1+\|z\|^2)(1+\|w\|^2)}
{|1+z\cdot w|^2}.
\]
Thus the radial holomorphic curve is $\gamma(z)=\bigl((z,0,\ldots,0),(0,\ldots,0)\bigr)$.
\item $M=\mathds H{\rm P}^n$. Using again Lee's affine coordinates $(z,w)$ on the Grassmannian realization of the complexification, take
\[
\gamma(z)=\bigl((z,0,\ldots,0),(0,\ldots,0)\bigr).
\]
In Lee's notation the invariant is $N=A/|B|^2$. Along this curve one has $A=1+|z|^2$, $B=-1$,
and therefore
\[
N(\gamma(z))=1+|z|^2.
\]

\item $M=\mathds O{\rm P}^2$. By \cite{Wolf}, $\mathds O{\rm P}^2$contains a totally geodesic copy of $\mathds H{\rm P}^2$. The adapted complex structure restricts to the complexification
of a totally geodesic submanifold (see \cite{lempert-szoke}). Moreover, the rank-one
radial coordinate $t$ restricts to the corresponding radial coordinate on
$\mathds H{\rm P}^2$. Hence the curve above for $\mathds H{\rm P}^2$, composed with the
holomorphic inclusion
\[
(\mathds H{\rm P}^2)^\C\hookrightarrow(\mathds O{\rm P}^2)^\C,
\]
gives a radial holomorphic curve in $(\mathds O{\rm P}^2)^\C$.
\end{enumerate}
\end{proof}
\begin{theorem}\label{thm:stenzel}
No positive multiple of the Stenzel Ricci-flat non-flat K\"ahler metric on the complexification of a compact rank-one symmetric space admits a local Kähler immersion into $\CP^N$, for any $N\leq\infty$.
\end{theorem}

\begin{proof}
We prove the statement for the simply connected compact
rank-one symmetric spaces. The real-projective case follows
from the Riemannian covering $S^d\to \R {\rm P}^d$, $d\geq2$.
Indeed, any local Kähler immersion of a positive multiple
of the Stenzel metric on the complexification of
${\R{\rm P}}^d$ into $\CP^N$ would induce one of its
pullback to the corresponding spherical complexification,
into the same target.

We apply Theorem \ref{thm:main1}. 
By Lemma \ref{curve}, Stenzel metric restricts to a metric on a radial holomorphic curve with potential given by
\begin{equation}\label{eq:stenzel-general-diastasis}
\Phi(x)=f(1+x),
\qquad x=|z|^2,
\end{equation}
and since from $N=\cosh^2t=1+x$ we get $x=\sinh^2 t$, it follows that 
$$
\rho(t)=f(\cosh^2t)=f(1+\sinh^2t)=\Phi(\sinh^2 t),
$$
and
\begin{equation}\label{giugno}
\rho'(t)=\sinh(2t)\Phi'(\sinh^2 t).
\end{equation}
By \cite[Example~2.4]{BD}, the Stenzel solution is even, smooth and strictly convex.
In particular, $\rho'(0)=0$ and $\rho'(t)>0$ for $t>0$. Integrating the
equation gives
\[
(\rho'(t))^d
=
dC\int_0^t\sinh^{p} u\,\sinh^{d-1-p}(2u)\,du.
\]
Near $t=0$ the right-hand side is $C_0t^d(1+O(t^2))$, with
$C_0>0$. Hence the regular branch satisfies
\[
\rho'(t)=AtH(t^2),
\qquad A>0,
\]
where $H$ is holomorphic near the origin and $H(0)=1$. By \eqref{giugno}, writing $\sinh(2t)=\sum_{j=1}^\infty \frac{(2t)^{2j-1}}{(2j-1)!}=2tK(t^2)$ with $K$ holomorphic and such that $K(0)\neq 0$, we get:
\begin{equation}\label{luglio}
\Phi'(\sinh^2t)=\frac{AtH(t^2)}{\sinh(2t)}=\frac{AH(t^2)}{2K(t^2)}.
\end{equation}
The right-hand side of \eqref{luglio} is an even holomorphic function of $t$. Moreover,
$$
x=\sinh^2t=t^2L(t^2), \qquad L(0)=1,
$$
so that $t^2$ is holomorphic in $x$ near the origin. Hence $\Phi'$,
and therefore $\Phi$, is holomorphic near $x=0$.
For $t>0$, analyticity follows from the analytic ODE and strict
convexity, since $x=\sinh^2t$ is a local real-analytic coordinate.

It remains to show that $\Phi'$ admits no single-valued meromorphic continuation to $\C$.

Set
\[
t_*:=\frac{i\pi}{2},
\qquad
x_*:=\sinh^2 t_*=-1,
\]
and define
\[
I(t):=dC\int_0^t \sinh^p u\,\sinh^q(2u)\,du,
\qquad
q=d-1-p.
\]
Thus
\[
\bigl(\rho'(t)\bigr)^d=I(t).
\]
We claim that the regular branch $\rho'=I^{1/d}$ extends
holomorphically along the segment joining $0$ to $t_*$. Indeed, for
$0<\theta\leq\pi/2$ we have
\[
I(i\theta)
=
dC\,i^d
\int_0^\theta
\sin^p s\,\sin^q(2s)\,ds
\neq 0.
\]
Here we used $p+q=d-1$, while the last inequality follows because the
integrand is positive for $0<s<\theta$. Therefore $I$ has no zeros
along the segment $i(0,\pi/2]$. Since the regular branch of $I^{1/d}$
is holomorphic at the origin, it admits a holomorphic continuation
along the whole segment and satisfies
\[
\rho'(t_*)\neq 0.
\]
Moreover, since $I(t_*)\neq 0$, this continuation is holomorphic in a
neighbourhood of $t_*$. Hence
\[
\rho'(t)=\rho'(t_*)+O(t-t_*).
\]
Furthermore, 
$$
\sinh(2t_*)=\sinh(i\pi)=0,
$$
 and
$$
\frac{d}{dt}\sinh(2t)|_{t=t_*}=2\cosh(2t_*)=2\cos(\pi)=-2,
$$
thus $\sinh(2t)=-2(t-t_*)+O((t-t_*)^3)$.  It follows that
$$
\Phi'(\sinh^2t)=\frac{\rho'(t_*)+O(t-t_*)}{-2(t-t_*)+O((t-t_*)^3)}=-\frac{\rho'(t_*)}{2(t-t_*)}+O(1).
$$
 Therefore, since 
 $$
 1+x=\cosh^2t=\cosh^2(t-t_*+t_*) =-\sinh^2 (t-t_*)=-(t-t_*)^2+O((t-t_*)^4),
 $$ 
 we get
 $$
 t-t_*=C(1+x)^{1/2}+O((1+x)^{3/2}), \quad C\neq 0,
 $$
 and
 \begin{equation}\label{stenzelphi}
 \Phi'(x)=-\frac{\rho'(t_*)}{2C}(1+x)^{-1/2}+O(1).
 \end{equation}
Hence $\Phi'(x)$ admits no single-valued meromorphic continuation to $\C$ (it has a branch point at $x_*$). Conclusion follows by Theorem \ref{thm:main1}.
\end{proof}

\subsection{Hwang--Singer scalar-flat ALE metrics}
Let $n\geq2$ and $k\geq1$.
The Hwang--Singer
momentum construction \cite{hwang-singer}, applied to the positive line
bundle $\mathcal O_{\CP^{n-1}}(k)$, gives a complete scalar-flat K\"ahler
metric on $\mathcal O_{\CP^{n-1}}(-k)$. 

Recall that the complement of the zero section of $\mathcal O_{\CP^{n-1}}(-k)$ is biholomorphic to $
(\C^n\setminus\{0\})/\mathbb Z_k$, where $\mathbb Z_k$ acts by scalar multiplication. The scalar-flat K\"ahler metric can be described on its covering $\C^n\setminus\{0\}$ as follows. 

Let $\tau$ be the momentum coordinate, normalized
so that the zero section corresponds to $\tau=1$. The Hwang--Singer
momentum profile is
\[
\varphi_{n,k}(\tau)
=
\frac{1}{k}
\left(
\tau+(k-n)\tau^{2-n}
+(n-k-1)\tau^{1-n}
\right),
\]
and satisfies
\[
\varphi_{n,k}(1)=0,
\qquad
\varphi_{n,k}'(1)=1.
\]
Writing $s=|z|^2$ on the covering, let $F_{n,k}(s)$ be a radial
K\"ahler potential for the lifted metric. The logarithmic fibre
coordinate is $k\log s$. After a positive homothety of the metric,
setting
\[
y(s):=sF'_{n,k}(s),
\]
we have $y=\tau$, and hence
\begin{equation}\label{518}
s y'(s)
=
k\varphi_{n,k}(y(s))
=
\frac{P_{n,k}(y(s))}{y(s)^{n-1}}.
\end{equation}where
\[
P_{n,k}(y)
=
y^n+(k-n)y+n-k-1.
\]
Moreover, the boundary conditions at the zero section give
\[
y(s)=1+C s^k+O(s^{2k}),
\]
for a constant $C>0$.
After a radial dilation we may normalize $C$,
without affecting the immersion problem.
Equivalently, on every branch on which $P_{n,k}(y)\neq0$,
\begin{equation}\label{519}
\frac{d\log s}{dy}
=
\frac{y^{n-1}}{P_{n,k}(y)}.
\end{equation}
Notice that for $n=2$ one has
\[
P_{2,k}(y)=(y-1)(y+k-1).
\]
Integrating \eqref{519}, and changing the radial coordinate by a
constant dilation, gives
\[
s^k=(y-1)(y+k-1)^{k-1}.
\]
In particular, $(n,k)=(2,1)$ gives the Burns--Simanca metric,
$(n,k)=(2,2)$ gives the Eguchi--Hanson metric, and the case $k=n$
gives Calabi's Ricci-flat metric on the canonical bundle
$\mathcal O_{\CP^{n-1}}(-n)$.
\begin{theorem}
Let $n\geq2$ and $k\geq1$, with $(n,k)\neq(2,1)$. No positive
multiple of the Hwang--Singer scalar-flat ALE metric on
$\mathcal O_{\CP^{n-1}}(-k)$ admits a local Kähler
immersion into $\CP^N$, for any $N\leq\infty$.
\end{theorem}
\begin{proof}
We apply Theorem \ref{thm:main1} to the following radial holomorphic curve. Fix $\lambda>0$ and consider, on the covering, the holomorphic curve
\[
\gamma(\zeta)=(\lambda,\zeta,0,\ldots,0).
\]
Any local Kähler immersion of the line-bundle metric would
pull back to an immersion of the covering metric into the
same target.
Setting $x:=|\zeta|^2$, a radial K\"ahler potential for the induced metric is given by:

\begin{equation}\label{520}
\Phi(x)=F_{n,k}\bigl(\lambda^2+x\bigr)
\end{equation}
On the
geometric branch $y>1$ we have
\[
P_{n,k}(1)=0, \qquad
P_{n,k}'(y)
=
n y^{n-1}+k-n
\geq k>0.
\]
It follows that $P_{n,k}(y)>0$ for every $y>1$. Moreover, as $y\to+\infty$,
\[
\frac{d\log s}{dy}
=
\frac{y^{n-1}}{P_{n,k}(y)}
=
\frac{1}{y}+O\left(\frac{1}{y^n}\right).
\]
Therefore for any fixed $y_0>1$
\[
\log\frac{s(y)}{s(y_0)}=\int^{y}_{y_0}\frac{t^{n-1}}{P_{n,k}(t)}\,dt=\log y+O(1) \qquad {\textrm as}\ y\rightarrow+\infty,
\]
and hence $s(y)\to+\infty$ as $y\to+\infty$. Thus the geometric
branch $y(s)$ is defined and analytic for every $s>0$. Consequently,
$F_{n,k}'(s)=y(s)/s$, and hence also $F_{n,k}(s)$, extends
analytically along the whole positive real axis. 

We now prove that $F'$ has no single-valued meromorphic
continuation to $\C$. Near $y=1$, equation \eqref{519} gives
\[
s=(y-1)^{1/k}h(y),
\qquad
h(1)\neq0,
\]
for a suitable holomorphic function $h$. Hence the inverse branch
$y(s)$ extends holomorphically through $s=0$ and satisfies
\[
y(s)=1+O(s^k).
\]
Starting from this branch, we analytically continue the inverse
relation along a path in the $y$-plane joining a point near $y=1$ to
$y=0$ and avoiding the remaining zeros of $P_{n,k}$. Let $s_*\neq0$
be the value of $s$ corresponding to $y=0$ along this continuation.
Assume first that $k\neq n-1$. Then
\[
P_{n,k}(0)=n-k-1\neq0.
\]
It follows from \eqref{519} that, as $y\to0$,
\[
\log\frac{s}{s_*}
=
\frac{y^n}{n(n-k-1)}
+
O(y^{n+1}).
\]
Consequently
\[
s-s_*
=
\frac{s_*}{n(n-k-1)}y^n
+
O(y^{n+1}),
\]
and the inverse branch satisfies
\[
y(s)
=
C(s-s_*)^{1/n}
+
O\bigl((s-s_*)^{2/n}\bigr),
\qquad
C\neq0.
\]
Since $n\geq2$, this gives non-trivial monodromy. As $s_*\neq0$, the
same non-integral power occurs in
\[
F_{n,k}'(s)=\frac{y(s)}s.
\]
It remains to consider $k=n-1$. In this case
\[
P_{n,n-1}(y)=y^n-y
=
y(y^{n-1}-1),
\]
and \eqref{519} becomes
\[
\frac{d\log s}{dy}
=
\frac{y^{n-2}}{y^{n-1}-1}.
\]
Hence, near $y=0$,
\[
\log\frac{s}{s_*}
=
-\frac{y^{n-1}}{n-1}
+
O(y^{2n-2}),
\]
and therefore
\[
y(s)
=
C(s-s_*)^{1/(n-1)}
+
O\bigl((s-s_*)^{2/(n-1)}\bigr),
\qquad
C\neq0.
\]
For $n\geq3$ this is again a genuine branch point. The remaining
case, that is $n=2$ and $k=1$, is the Burns--Simanca metric and is
excluded from the statement. 
Thus, for every $(n,k)\neq(2,1)$, the radial potential extends
analytically along the positive real axis whereas its derivative has
no single-valued meromorphic continuation to $\C$. Conclusion follows by Theorem \ref{thm:main1}.
\end{proof}

\subsection{Hwang--Singer scalar-flat metrics on $\mathbb C^{n+1}$}

Let $n\geq1$ and $\beta<0$. The Hwang--Singer momentum
construction \cite{hwang-singer}, applied to the trivial line bundle
over $\mathbb C^n$ endowed with its flat metric, gives a
complete scalar-flat Kähler metric $\omega_\beta$ on
$\mathbb C^{n+1}$. In coordinates
$(z,\xi)\in\mathbb C^n\times\mathbb C$, a Kähler potential
away from $\xi=0$ is
\[
K_\beta(z,\xi)
=
\|z\|^2+
4f\left(\frac12\log\left(
|\xi|^2e^{-\frac{\beta}{2}\|z\|^2}
\right)\right),
\]
where
\[
\tau=f'(t),
\qquad
f''(t)=\varphi(\tau),
\qquad
\varphi(\tau)=\frac{2\tau}{(1-\beta\tau)^n}.
\]
The metric extends smoothly across $\xi=0$.
In \cite{cristofori-zedda}, it was proved that no positive multiple of
$\omega_\beta$ is projectively induced.

\begin{theorem}
Let $n\geq1$ and $\beta<0$.
No positive multiple of the Hwang--Singer scalar-flat
metric $\omega_\beta$ on $\mathbb C^{n+1}$ admits a local
Kähler immersion into $\CP^N$, for any $N\leq\infty$.
\end{theorem}

\begin{proof}
We apply Theorem~\ref{thm:main1} to the holomorphic curve
\[
\gamma(\zeta)=(0,\ldots,0,\zeta).
\]
Setting $x:=|\zeta|^2$, the restricted potential is
\[
\Phi(x)=4f\left(\frac12\log x\right),
\qquad x>0.
\]
Since $t=\frac12\log x$ and $\tau=f'(t)$, we have
\[
x\Phi'(x)=2\tau(x),
\qquad
\frac{d\tau}{dx}
=
\frac{\tau}{x(1-\beta\tau)^n}.
\]
Equivalently,
\[
\frac{d\log x}{d\tau}
=
\frac{(1-\beta\tau)^n}{\tau}.
\]
Integrating and applying a constant dilation of the fibre
coordinate, we obtain
\[
x=\tau e^{P(\tau)},
\qquad
P(\tau)=
\sum_{j=1}^{n}\binom{n}{j}\frac{(-\beta)^j}{j}\tau^j,
\]
and differentiating we get
\[
\frac{dx}{d\tau}
=
e^{P(\tau)}(1-\beta\tau)^n.
\]
For $\tau\geq0$, this derivative is positive and $x(\tau)$
increases from $0$ to $+\infty$. Hence its inverse $\tau(x)$
is analytic along the whole non-negative real axis.
Moreover, $x=\tau+O(\tau^2)$ near $\tau=0$, so
\[
\tau(x)=x+O(x^2).
\]
It follows that $\Phi'(x)=2\tau(x)/x$ extends holomorphically
through $x=0$. Thus $\Phi$ also extends holomorphically
through $x=0$ and along the whole positive real axis.

It remains to show that $\Phi'$ has no single-valued
meromorphic continuation to $\mathbb C$.
The function $x(\tau)=\tau e^{P(\tau)}$ is entire and its
derivative vanishes only at
\[
\tau_*=\frac1\beta<0,
\]
where it has a zero of order $n$. The inverse branch
$\tau(x)$ can be continued from $\tau=0$ towards $\tau_*$
along the negative real axis without meeting another
critical point. Set
\[
x_*=\tau_*e^{P(\tau_*)}\neq0.
\]
Near $\tau_*$ we have
\[
x-x_*
=
A(\tau-\tau_*)^{n+1}
+
O\bigl((\tau-\tau_*)^{n+2}\bigr),
\qquad A\neq0,
\]
and therefore
\[
\tau(x)
=
\tau_*+C(x-x_*)^{1/(n+1)}
+
O\bigl((x-x_*)^{2/(n+1)}\bigr),
\qquad C\neq0.
\]
Since $x_*\neq0$, the same non-integral power occurs with
nonzero coefficient in
\[
\Phi'(x)=\frac{2\tau(x)}{x}.
\]
Thus $\Phi'$ has non-trivial monodromy around $x_*$
and cannot admit a single-valued meromorphic continuation
to $\mathbb C$. The conclusion follows from Theorem~\ref{thm:main1}.
\end{proof}

\section{Invariant rational curves and moment maps}\label{sec:moment}
In this section we prove Theorem~\ref{thm3}. We first establish
a necessary condition for an invariant metric on
$\mathbb{CP}^1$ to be projectively induced, expressed
in terms of its moment map.

\begin{lemma}\label{cp1}
Let $S^1$ act non-trivially and holomorphically on $\CP^1$  and let $\omega$ be a $S^1$-invariant K\"ahler form. If $\omega$ is induced by a K\"ahler immersion into $\CP^N$, $N\leq \infty$, then the moment map for the action is rational.
\end{lemma}
\begin{proof}
Since $\omega$ is $S^1$-invariant and projectively induced, the standard
description of projectively induced toric metrics (see, e.g.,
\cite[Lemma~4.1]{ALZ}) gives, in a suitable affine coordinate,
\[
\omega=\frac{i}{2}\partial\bar\partial\log P(|z|^2),
\qquad
P(\rho)=\sum_{k=0}^{\ell}a_k\rho^k,
\]
where $a_k\geq0$ and $a_0a_\ell>0$.
In particular, using that $\partial\bar \partial h(\rho)=(h'(\rho)+\rho h''(\rho))dz\wedge d\bar z$, we get
$$
\omega=\frac i2 (\rho h'(\rho))' dz\wedge d\bar z,
$$
with $h(\rho)=\log P(\rho)$.

Let $X=im\left(z\partial_z-\bar z \partial_{ \bar z}\right)$ be the infinitesimal generator of the action and denote by $\mu$ the moment map, defined by $d\mu=-\iota_X\omega$.
Then
$$
\iota_Xdz\wedge d\bar z=mi (zd\bar z+\bar zdz)=mid\rho.
$$
and thus
$$
d\mu=-\iota_X\omega=\frac m2 (\rho h'(\rho))'d\rho
$$
which implies
$$
\mu=\frac m2 \rho h'(\rho)+c=\frac m2 \frac{\sum_{k=0}^\ell ka_k \rho^{k}}{\sum_{k=0}^\ell a_k \rho^{k}}+c
$$
for some constant $c$. 

Thus \(\mu\) is a rational function of \(\rho=|z|^2\) in an affine coordinate adapted to the action, as claimed.
\end{proof}

\begin{proof}[Proof of Theorem~\ref{thm3}]
Assume that $c\omega$ is projectively induced for some $c>0$. Then so is its restriction to $\Sigma$ and, since $\Sigma\simeq \CP^1$ is simply connected, the immersion is global. By Lemma~\ref{cp1}, its moment map $c\mu$ is rational. Since $c>0$, also $\mu$ is rational, contradicting the hypothesis.
\end{proof}

We conclude this section with the following proposition, that relates the non-rationality of the moment map to the analytic obstruction in Theorem \ref{thm:main1}, and will be used in Section \ref{7}.
\begin{proposition}\label{applyalgebraic}
Let $(M,\omega)$ be a real-analytic K\"ahler manifold containing a holomorphic curve $\Sigma\simeq{\CP}^1$ whose induced metric is invariant under a non-trivial holomorphic
$S^1$-action. Then, at each fixed point of this action, the corresponding radial potential $\Phi$ extends analytically along $[0,+\infty)$.
If moreover the moment map $\mu$, viewed as a function of $\rho=|z|^2$ in an affine coordinate adapted to the action, is non-rational, then $\Phi'$ has no single-valued meromorphic continuation to $\C$.
\end{proposition}
\begin{proof}
A non-trivial holomorphic $S^1$-action on $\Sigma\simeq\mathbb{CP}^1$ has exactly two fixed points. Choose an affine coordinate $z$ in which these points
are $0$ and $\infty$ and the action is $z\mapsto e^{im\theta}z$,  $m\in\Z\setminus\{0\}$.
Set $\rho=|z|^2$.

Let $\Phi(|z|^2)$ be the radial potential near $z=0$ of the induced metric on $\Sigma$. Then as in the proof of Lemma \ref{cp1}, we have:
$$
\mu(\rho)=\frac m2\rho\,\Phi'(\rho)+\mu(0)
$$
and thus
\begin{equation}\label{rhoalg}
\Phi(\rho)=\frac2m\int_0^\rho\frac{\mu(t)-\mu(0)}{t}dt+\Phi(0).
\end{equation}
Since $\mu(\rho)$ is the moment map of a real analytic $S^1$-invariant metric on $\CP^1$, it is $S^1$ invariant on the affine chart $\C$ and hence it is a real analytic function of $\rho$. Therefore, for every $\rho_0\geq 0$, $\mu$ extends holomorphically to a complex neighborhood of $\rho_0$. Moreover, $\mu(\rho)-\mu(0)=O(\rho)$ near $0$ and hence by \eqref{rhoalg} also $\Phi(\rho)$ extends holomorphically through every point of the positive real axis.

To prove the last assertion, suppose that $\Phi'$ admits a single-valued meromorphic continuation to $\C$. Then the identity
$$
\mu(\rho)=\frac{m}{2}\rho\Phi'(\rho)+\mu(0)
$$
gives such a continuation for $\mu$.
Consider now the other fixed point of the action on $\Sigma$.
In the coordinate $w=1/z$, the moment map is real analytic
and rotation invariant. Hence it has a convergent expansion
\[
\mu=\mu_\infty+\sum_{j\geq1}a_j|w|^{2j},
\]
and therefore
\[
\mu(\rho)=\mu_\infty+\sum_{j\geq1}a_j\rho^{-j}
\]
for sufficiently large positive $\rho$.
This series defines a holomorphic function in a neighbourhood
of infinity. By the identity theorem, it agrees with the
meromorphic continuation of $\mu$ on their common domain.
Thus $\mu$ extends meromorphically to $\mathbb{CP}^1$
and is therefore rational, contradicting the hypothesis.
\end{proof}

\section{Gibbons--Hawking ALE metrics}\label{sec:GH}
For $k\geq1$, let $a_0,\ldots,a_k\in\R^3$ be distinct points, called the centres, and consider the harmonic function
$$
V\!:\R^3\setminus\{a_0,\ldots,a_k\}\longrightarrow\R_{>0},\qquad V(q)=\frac12\sum_{j=0}^k\frac1{|q-a_j|}.
$$
The Gibbons--Hawking construction associates with $V$ a complete ALE hyper-K\"ahler metric on a four-dimensional manifold $M$, given away from the fixed points of the circle action by
$$
g=V\,g_{\R^3}+\frac1V\eta^2,
$$
where $\eta$ is a connection one-form on the corresponding principal $S^1$-bundle over $\R^3\setminus\{a_0,\ldots,a_k\}$. Observe that away from the fixed points, $M$ is an $S^1$-fibration over $\R^3\setminus\{a_0,\ldots,a_k\}$, and the circle action is tri-holomorphic; its fixed points project precisely to the centres $a_0,\ldots,a_k$.

The following argument begins with an arbitrary direction $u$ and proves the required obstruction in $I_u$ directly. In particular, it treats collinear and non-collinear configurations, including two centres, in a single proof. Theorem~\ref{thm:main4}(vi) then follows from Theorem~\ref{thm:main1}.

\begin{lemma}\label{algebraicGH}
Let $a_0,\ldots,a_k\in\R^3$ be distinct, with $k\geq1$, and fix any $u\in S^2$. For the compatible complex structure $I_u=u_1I+u_2J+u_3K$, the corresponding Gibbons--Hawking manifold contains an invariant holomorphic curve through a circle fixed point, biholomorphic to $\C$, with a globally defined radial coordinate $w$. Its radial K\"ahler potential $\Phi(\rho)$, $\rho=|w|^2$, extends real-analytically along $[0,\infty)$, and $\Phi'$ is algebraic and non-rational. Consequently $\Phi'$ has no single-valued meromorphic continuation to $\C$.
\end{lemma}
\begin{proof}
Choose $a_0$ among the centres so that $\langle a_0,u\rangle$ is maximal, and put $q(x)=a_0+xu$ for $x\geq0$. If another centre $a_j$ lay on this ray, then $a_j=a_0+xu$ with $x>0$ and $\langle a_j,u\rangle>\langle a_0,u\rangle$, a contradiction. Thus the open ray avoids every centre, even when several centres share the same maximal projection. Let $\mu:M\to\R^3$ be the Gibbons--Hawking projection, and let $C_u$ be the preimage of the open ray together with the fixed point over $a_0$. The Gibbons--Hawking expression for $I_u$ sends the horizontal lift of $u$ to the vertical circle generator, up to positive metric factors and the consistent choice of orientation (see, e.g., \cite{LotayOliveira}). Hence the tangent plane of $C_u$ is $I_u$-invariant. The standard smooth extension at a centre makes $C_u$ an embedded $I_u$-holomorphic curve at its fixed point as well.

Trivialize the circle bundle along the open ray so that $\eta=d\tau$, where $\tau$ has period $2\pi$. The induced metric and K\"ahler form are
$$
g_{C_u}=V(q(x))\,dx^2+V(q(x))^{-1}d\tau^2,\qquad \omega_{C_u}=dx\wedge d\tau.
$$
Define $s$ by $ds=V(q(x))\,dx$ and set $w=e^{s+i\tau}$, $\rho=|w|^2=e^{2s}$. The sign of $\tau$ is chosen so that $w$ is holomorphic. Then
\begin{equation}\label{eq:GHoutward}
\frac{d\log\rho}{dx}=2V(q(x)),\qquad \rho'(x)=2V(q(x))\rho(x)>0.
\end{equation}
Near $x=0$ we have $2V(q(x))=x^{-1}+h(x)$, with $h$ real analytic. Consequently $\rho(x)=C_0x(1+O(x))$, $C_0>0$. The smooth fixed-point model for the circle action shows that $w$ extends across the added point as a holomorphic coordinate, with $w=0$ there. We will see below that $\rho(x)\to\infty$ as $x\to\infty$, so $w$ identifies $C_u$ biholomorphically with $\C$.

For $j\geq1$ write $a_j-a_0=c_ju+v_j$, where $c_j=\langle a_j-a_0,u\rangle\leq0$, $v_j\perp u$. 
Distinctness of the centres implies $c_j<0$ whenever $|v_j|=0$. Therefore
$$
2V(q(x))=\frac1x +\sum_{\substack{1\leq j\leq k\\ |v_j|=0}}\frac1{x-c_j}+\sum_{\substack{1\leq j\leq k\\ |v_j|>0}}\frac1{\sqrt{(x-c_j)^2+|v_j|^2}}.
$$
Integrating \eqref{eq:GHoutward} gives, for some $C>0$,
\begin{equation}\label{eq:GHoutwardrho}
\rho(x)=Cx\prod_{\substack{1\leq j\leq k\\|v_j|=0}}(x-c_j)\prod_{\substack{1\leq j\leq k\\|v_j|>0}}\left(x-c_j+\sqrt{(x-c_j)^2+|v_j|^2}\right).
\end{equation}
All factors in the products are positive and real analytic for $x\geq0$. Formula \eqref{eq:GHoutwardrho} makes $\rho(x)$ algebraic and yields $\rho(x)\sim C_\infty x^{k+1}$ as $x\to\infty$ for some $C_\infty>0$. Together with \eqref{eq:GHoutward}, this proves that $\rho:[0,\infty)\to[0,\infty)$ is strictly increasing and has a real-analytic inverse $x=x(\rho)$ throughout the half-line. The inverse is algebraic, since exchanging the two variables in a polynomial relation for $\rho(x)$ gives one for $x(\rho)$. But
$$
x(\rho)\sim (\rho/C_\infty)^{1/(k+1)}\qquad(\rho\to+\infty),
$$
and $1/(k+1)$ is non-integral for $k\geq1$. A rational function has integral order of growth at infinity, so $x(\rho)$ is non-rational.

A radial potential for $\omega_{C_u}=\frac{i}{2}\partial\bar\partial\Phi(\rho)$ is
$$
\Phi(\rho)=2\int_0^\rho\frac{x(t)}{t}\,dt,\qquad \Phi'(\rho)=\frac{2x(\rho)}{\rho}.
$$
Indeed, $x(0)=0$, so the integrand is regular at the origin, and $\rho\Phi'(\rho)=2x(\rho)$ gives $\omega_{C_u}=dx\wedge d\tau$. The inverse-function argument above makes $\Phi$ real analytic along the positive axis. Its derivative is algebraic and non-rational: if $\Phi'$ were rational, then $x(\rho)=\rho\Phi'(\rho)/2$ would be rational. A single-valued meromorphic algebraic function on $\C$ is rational, so $\Phi'$ has no such continuation. Theorem~\ref{thm:main1} thus applies to $(M,\omega_{I_u})$ and rules out a local K\"ahler immersion of every positive multiple into every $\CP^N$, $N\leq\infty$. Since $u$ was arbitrary from the start, the conclusion holds for all compatible complex structures in one step.
\end{proof}

\section{K\"ahler immersions into other complex space forms}
\label{7}

We collect here the extensions of our results to complex space
forms of arbitrary curvature and signature.

For $b\in\mathbb{R}$ and
$r,s\in\mathbb{Z}_{\geq 0}\cup\{\infty\}$, denote by $X_b^{r,s}$
the complex space form of constant holomorphic sectional curvature
$4b$ and complex signature $(r,s)$, where $r$ and $s$ count the
positive and negative directions, respectively. We refer to
\cite{BR,R,Uind} for indefinite complex space forms and the
related Calabi theory. In particular,
\[
X_{-1}^{N,0}=\mathbb{CH}^N,\qquad
X_0^{N,0}=\mathbb{C}^N,\qquad
X_1^{N,0}=\mathbb{CP}^N.
\]
The diastasis at the origin is
\[
D_b(u,v)=
\begin{cases}
\displaystyle
\frac{1}{b}\log\left(
1+b\left(\sum_{j=1}^{r}|u_j|^2
-\sum_{k=1}^{s}|v_k|^2\right)\right),
& b\neq 0,\\[8pt]
\displaystyle
\sum_{j=1}^{r}|u_j|^2-\sum_{k=1}^{s}|v_k|^2,
& b=0.
\end{cases}
\]
Accordingly, set
\[
E_b(f)=
\begin{cases}
\dfrac{e^{bf}-1}{b},& b\neq 0,\\[6pt]
f,& b=0.
\end{cases}
\]

In the positive-definite setting, Bochner's constructions
\cite{bochner} give K\"ahler immersions
\[
\mathbb{CH}^{N}\longrightarrow
\mathbb{C}^{\infty}\longrightarrow
\mathbb{CP}^{\infty}.
\]
Thus a K\"ahler immersion into a hyperbolic space form yields
one into a flat space form, and a K\"ahler immersion into a flat
space form yields one into projective space. Since our
non-immersion results hold for every positive rescaling, they
therefore exclude all positive-definite complex space forms,
of any curvature and finite or infinite dimension.

The indefinite setting has a different feature: allowing both
signature indices to be infinite gives universal local targets.

\begin{proposition}
\label{prop:universal-indefinite}
Every real-analytic K\"ahler manifold admits a local K\"ahler
immersion into $X_b^{\infty,\infty}$ for every
$b\in\mathbb{R}$.
\end{proposition}

\begin{proof}
Calabi's local immersion theorem \cite{calabi} gives
\[
D_p=\sum_j|f_j|^2-\sum_k|g_k|^2,
\]
where the holomorphic functions vanish at $p$ and both sums
converge locally. This proves the assertion for $b=0$.

For $b\neq 0$, expand $E_b(D_p)$ and separate the positive and
negative terms. This gives holomorphic functions $F_\alpha$
and $G_\beta$, vanishing at $p$, such that
\[
E_b(D_p)=
\sum_\alpha|F_\alpha|^2-\sum_\beta|G_\beta|^2.
\]
For $b>0$, the signs are determined by the parity of the total
degree in the functions $g_k$; for $b<0$, positive terms
correspond to odd total degree in the functions $f_j$, and
negative terms to even total degree. Both families converge
locally, as follows from the absolute convergence of the
exponential series with argument
\[
|b|\left(\sum_j|f_j|^2+\sum_k|g_k|^2\right).
\]
Consequently,
\[
D_p=\frac{1}{b}\log\left(
1+b\left(
\sum_\alpha|F_\alpha|^2-\sum_\beta|G_\beta|^2
\right)\right),
\]
which gives the required local K\"ahler immersion.
\end{proof}

Consequently, no real-analytic K\"ahler metric can be excluded
from all indefinite complex space forms. The relevant question
is whether an immersion exists with at least one signature
index finite. Unlike the positive-definite situation, passing
between target geometries does not in general preserve this
condition: the expansion above may produce infinitely many
terms of both signs even when the initial decomposition has
only finitely many terms of one sign. Thus projective
non-immersion alone does not exclude indefinite flat or
hyperbolic immersions with a finite index. These targets are
instead treated directly by the same Pringsheim argument.

\begin{proposition}
\label{prop:finite-index-obstruction}
Under the hypotheses of Theorem~\ref{thm:main1}, for every
$c>0$ and every $b\in\mathbb{R}$, the K\"ahler manifold
$(M,c\omega)$ admits no local K\"ahler immersion into
$X_b^{r,s}$ unless $r=s=\infty$.
\end{proposition}

\begin{proof}
Replace the exponential in the proof of
Theorem~\ref{thm:main1} by $E_b(c\Phi)$.
The coefficient comparison in Lemma~\ref{lem:diagonal}
applies equally to this function. An immersion into
$X_b^{r,s}$ forces its diagonal coefficient matrix on the
radial curve to have at most $r$ positive and at most $s$
negative entries. This follows from the representation of
$E_b(cD_p)$ as a sum of $r$ squared moduli minus a sum of $s$
squared moduli: every finite principal submatrix satisfies
the corresponding signature bounds.

The logarithmic-derivative argument used in
Theorem~\ref{thm:main1} shows that $E_b(c\Phi)$ cannot be
entire when $b\neq 0$, since
$1+bE_b(c\Phi)=e^{bc\Phi}$. For $b=0$, entire continuation
of $E_0(c\Phi)=c\Phi$ would make $\Phi'$ entire.
Thus $E_b(c\Phi)$ has finite convergence radius and continues
analytically along the whole positive real axis.

Pringsheim's Theorem, after subtracting finitely many Taylor
terms, excludes eventual non-negativity of its coefficients.
Applying the same argument to $-E_b(c\Phi)$ excludes eventual
non-positivity. Hence infinitely many positive-degree
coefficients have each sign, and both $r$ and $s$ must be
infinite.
\end{proof}

Combining the preceding propositions gives the following
extensions of our main results.

\begin{corollary}
\label{cor:all-space-forms}
Let $\omega$ be any metric listed in
Theorem~\ref{thm:main4}, or let $(M,\omega)$ satisfy the
hypotheses of Theorem~\ref{thm3}. For every $c>0$ and every
$b\in\mathbb{R}$, the metric $c\omega$ admits a local
K\"ahler immersion into $X_b^{r,s}$ if and only if
$r=s=\infty$.

For the Gibbons--Hawking metrics, this statement holds in
every compatible complex structure.
\end{corollary}

\begin{proof}
The proofs of the applications already verify the analytic
hypotheses of Theorem~\ref{thm:main1}, and hence of
Proposition~\ref{prop:finite-index-obstruction}.
For the invariant rational curves of Theorem~\ref{thm3},
the same hypotheses follow from
Proposition~\ref{applyalgebraic}.
Existence when both indices are infinite follows from
Proposition~\ref{prop:universal-indefinite}.
\end{proof}

\subsection*{Acknowledgements}
The authors are members of GNSAGA of INdAM.

\subsection*{AI Disclosure} During the preparation of this work, the authors used {\it Claude Opus 4.8} and {\it ChatGPT-6 Astra} for language editing and review. All the mathematical results are due to the authors.



\begin{thebibliography}{99}

\bibitem{ALL}
C. Arezzo, C. Li and A. Loi,
\emph{Gromov-Hausdorff limits and holomorphic isometries}, {\em Math. Res. Lett.} {\bf 32} (2025), Number 3, 707-738.

\bibitem{ALPZradial} C. Arezzo, A. Loi, G. Placini and M. Zedda, \emph{Radial projectively induced canonical K\"ahler metrics: rigidity and classification}, arXiv:2609.08438 (2026).

\bibitem{ALZ} C. Arezzo, A. Loi, F. Zuddas, On homothetic balanced metrics, {\em Ann. Global Anal. Geom.} {\bf41} (2012), 473--491.

\bibitem{BR} M. Barros, A. Romero, Indefinite K\"ahler manifolds, {\em Math. Ann.} {\bf261} (1982), 55--62.


\bibitem{BD} O. Biquard, T. Delcroix, Ricci flat K\"ahler metrics on rank two complex symmetric spaces, {\em J. \'Ec. polytech. Math.} {\bf6} (2019), 163--201.

\bibitem{bochner} S. Bochner, Curvature in Hermitian metric, {\em Bulletin of the Am. Math. Soc.} {\bf53} (1947), 179--195.

\bibitem{calabi} E. Calabi, \emph{Isometric imbedding of complex manifolds}, Ann. of Math.
(2) \textbf{58} (1953), 1--23.

\bibitem{cdlo} P. Candelas and X. C. de la Ossa, \emph{Comments on conifolds}, Nuclear Phys.
B \textbf{342} (1990), 246--268.

\bibitem{cannas-loi} F. Cannas Aghedu and A. Loi, \emph{The Simanca metric admits a regular
quantization}, Ann. Global Anal. Geom. \textbf{56} (2019), no. 3, 583--596.


\bibitem{cristofori-zedda} S. Cristofori and M. Zedda, \emph{K\"ahler geometry of scalar flat
metrics on line bundles over polarized K\"ahler--Einstein manifolds}, J. Geom. Anal.
\textbf{34} (2024), doi:10.1007/s12220-024-01590-0.

\bibitem{flajolet-sedgewick} P. Flajolet and R. Sedgewick, \emph{Analytic Combinatorics}, Cambridge University Press, 2009. Theorem IV.6,

\bibitem{DSIL} A. J. Di Scala, H. Ishi,  A. Loi, K\"ahler Immersions of Homogeneous K\"ahler Manifolds into Complex Space Forms, {\em Asian J. Math.} {\bf16} (2012), no. 3, 479--488.

\bibitem{D} I. Dolgachev, {\em Lectures on Invariant Theory}, London Mathematical Society Lecture Note Series 296, Cambridge University Press, 2003.

%


\bibitem{GH} G. W. Gibbons and S. W. Hawking, Gravitational multi-instantons, {\em Phys. Lett. B} {\bf78} (1978), 430--432.

\bibitem{HuangTaubNut} S. Huang, \emph{The Taub--NUT metric is not projectively induced}, arXiv:2609.21950 (2026).

\bibitem{hwang-singer} A. D. Hwang and M. A. Singer, \emph{A momentum construction for circle-invariant K\"ahler metrics}, Trans. Amer. Math. Soc. \textbf{354} (2002), no.~6, 2285--2325.

\bibitem{Hulin}
D.~Hulin,
\emph{Kähler–Einstein metrics and projective embeddings},
J. Geom. Anal. \textbf{10} (2000), 525–528.





\bibitem{lebrun} C. LeBrun, \emph{Complete Ricci-flat K\"ahler metrics on $\C^n$ need not be
flat}, in Several Complex Variables and Complex Geometry, Part 2, Proc. Sympos. Pure Math.
\textbf{52}, Amer. Math. Soc., Providence, RI, 1991, 297--304.


\bibitem{Lee}
T.-C. Lee,
\emph{Complete Ricci flat K\"ahler metric on $M_I^n$, $M_{II}^{2n}$, $M_{III}^{4n}$},
Pacific J. Math. \textbf{185} (1998), no.~2, 315--326.

\bibitem{lempert-szoke}
L. Lempert and R. Sz\"oke,
\emph{Global solutions of the homogeneous complex Monge--Amp\`ere equation
and complex structures on the tangent bundle of Riemannian manifolds},
Math. Ann. \textbf{290} (1991), 689--712.


\bibitem{lsz} A. Loi, F. Salis, and F. Zuddas, \emph{Two conjectures on Ricci-flat K\"ahler
metrics}, Math. Z. \textbf{290} (2018), 599--613.

\bibitem{lz-book} A. Loi, M. Zedda, \emph{K\"ahler Immersions of K\"ahler Manifolds into
Complex Space Forms}, Lecture Notes of the Unione Matematica Italiana \textbf{23}, Springer,
Cham, 2018.

\bibitem{loi-zedda-wallach} A. Loi, M. Zedda, \emph{K\"ahler--Einstein submanifolds of the
infinite dimensional projective space}, Math. Ann. \textbf{350} (2011), no. 1, 145--154.

\bibitem{lzz-taubnut} A. Loi, M. Zedda, and F. Zuddas, \emph{Some remarks on the K\"ahler
geometry of the Taub--NUT metrics}, Ann. Global Anal. Geom. \textbf{41} (2012), 515--533.

\bibitem{lzz-calabi} A. Loi, M. Zedda, and F. Zuddas, \emph{Ricci flat Calabi's metric is not
projectively induced}, Tohoku Math. J. (2) \textbf{73} (2021), no. 1, 29--37.

\bibitem{LotayOliveira}
J. D. Lotay, G. Oliveira, {Special Lagrangians, Lagrangian mean curvature flow and the Gibbons--Hawking ansatz}, {\em J. Diff. Geom.} {\bf126} (2024), no.~3, 1121--1184.

\bibitem{MTZ} S. Marini, N. Tardini, M. Zedda, Projectively Induced K\"ahler Cones over Regular Sasakian Manifolds, {\em Geom. Dedicata} {\bf218} (2024), article 85.

\bibitem{R} A. Romero, An extension of Calabi's rigidity theorem to complex submanifolds of indefinite complex space forms, {\em Manuscripta Math.} {\bf59} (1987), 261--276.


\bibitem{stenzel} M. B. Stenzel, \emph{Ricci-flat metrics on the complexification of a
compact rank one symmetric space}, Manuscripta Math. \textbf{80} (1993), 151--163.

\bibitem{Uind} M. Umehara, {\em Diastases and real analytic functions on complex manifolds}, J. Math. Soc. Japan {\bf40} (1988), 519--539.


\bibitem{TW} S. Tolman, J. Watts, Tame circle actions, {\em Trans. of the AMS} {\bf 369} (2017), 7443--7467.

\bibitem{Wolf}
J.~A. Wolf,
Elliptic spaces in Grassmann manifolds, {\em Illinois J. Math.} {\bf7} (1963), 447--462.


\bibitem{zedda-stenzel} M. Zedda, \emph{Stenzel's Ricci-flat K\"ahler metrics are not
projectively induced}, Osaka J. Math. \textbf{58} (2021), no. 4, 921--927.

\end{thebibliography}
\end{document}